\documentclass[12pt]{article}

\usepackage{amsmath}
\usepackage{amsthm}
\usepackage{graphics}
 
\newtheorem{theorem}{Theorem}
\newtheorem{cor}[theorem]{Corollary}
\newtheorem*{pf1}{Proof of Theorem~\ref{main}}
\newtheorem*{pf2}{Proof of Corollary~\ref{maincor}}
\newtheorem*{pf3}{Proof of Theorem~\ref{theorem}}
\newtheorem*{ack}{Acknowledgement}

\begin{document}

\bibliographystyle{amsplain}

\begin{center}
\LARGE{A Four-parameter Partition Identity} 
\bigskip

\large Cilanne E. Boulet \\
\small Department of Mathematics \\
Massachusetts Institute of Technology\\
Cambridge, MA, 02139 \\
\emph{e-mail:} cilanne@math.mit.edu \\
\end{center}

\normalsize
\section{Introduction}
\label{sec-intro}
In \cite{Andr}, Andrews considers partitions with respect to
size, number of odd parts, and number of odd parts of the conjugate.  
He derives the following generating function
\begin{equation}
\label{And:result}
\sum_{\lambda \in \mathrm{Par}} r^{\theta(\lambda)} s^{\theta(\lambda')} 
q^{|\lambda|} = 
\prod_{j=1}^ \infty 
\frac{(1+rsq^{2j-1})}{(1-q^{4j})(1-r^2q^{4j-2})(1-s^2q^{4j-2})}
\end{equation}
where $\mathrm{Par}$ denotes the set of all partitions, 
$|\lambda|$ denotes the
size (sum of the parts) of $\lambda$, $\theta(\lambda)$ denotes the
number of odd parts in the partition $\lambda$, and $\theta(\lambda')$ denotes 
the number of odd parts in the conjugate of $\lambda$.  
A combinatorial proof of Andrews' result was found by Sills in \cite{Sills}.  

In this paper, we generalize this result and provide a combinatorial proof
of our generalization.  This gives a simpler combinatorial proof
of (\ref{And:result}) than the one found in \cite{Sills}.  

\section{Main Result}
\label{sec-main}
Let $\lambda = (\lambda_1, \lambda_2, \dots )$ be a partition of $n$, denoted
$\lambda \vdash n$.  Consider the following weight functions on the set of all
partitions:
\begin{eqnarray*}
\alpha(\lambda) & = & \sum \lceil \lambda_{2i-1}/2 \rceil  \\   
\beta(\lambda)  & = & \sum \lfloor \lambda_{2i-1}/2 \rfloor \\
\gamma(\lambda) & = & \sum \lceil \lambda_{2i}/2 \rceil  \\
\delta(\lambda) & = & \sum \lfloor \lambda_{2i}/2 \rfloor .
\end{eqnarray*}
Also, let $a$, $b$, $c$, $d$ be (commuting) indeterminants, and define
\[
w(\lambda) = a^{\alpha(\lambda)}b^{\beta(\lambda)}c^{\gamma(\lambda)}
d^{\delta(\lambda)}.
\]
For instance, if $\lambda = (5,4,4,3,2)$ then $\alpha(\lambda)$ is the number of
$a$'s in the following diagram for $\lambda$, $\beta(\lambda)$ is the number of
$b$'s in the diagram, $\gamma(\lambda)$ is the number of
$c$'s in the diagram, and $\delta(\lambda)$ is the number of
$d$'s in the diagram.  Moreover, $w(\lambda)$ is the product of the
entries of the diagram.
\[
\begin{array}{ccccc}
a & b & a & b & a \\ 
c & d & c & d &   \\
a & b & a & b &   \\
c & d & c &   &   \\
a & b &   &   &   \\
\end{array}
\]
These weights were first suggested by Stanley in \cite{Stan}.

Let $\Phi(a,b,c,d) = \sum w(\lambda)$, where the sum is over all
partitions $\lambda$, and let $\Psi(a,b,c,d) = \sum w(\lambda)$, 
where the sum is over all
partitions $\lambda$ with distinct parts.  
We obtain the following product formulas for 
$\Phi(a,b,c,d)$ and $\Psi(a,b,c,d)$:  

\begin{theorem}
\label{main}
\[
\Phi(a,b,c,d) = \prod_{j=1}^\infty 
\frac{(1+a^{j}b^{j-1}c^{j-1}d^{j-1})(1+a^{j}b^{j}c^{j}d^{j-1})}
{(1-a^{j}b^{j}c^{j}d^{j})(1-a^{j}b^{j}c^{j-1}d^{j-1})
(1-a^{j}b^{j-1}c^{j}d^{j-1})}
\]
\end{theorem}

\begin{cor}
\label{maincor}
\[
\Psi(a,b,c,d) = \prod_{j=1}^\infty 
\frac{(1+a^{j}b^{j-1}c^{j-1}d^{j-1})(1+a^{j}b^{j}c^{j}d^{j-1})}
{(1-a^{j}b^{j}c^{j-1}d^{j-1})}
\]
\end{cor}

Andrews' result follows easily from Theorem~\ref{main}.  
Note that we can express number of odd parts of $\lambda$, number of odd
parts of $\lambda'$ and size of $\lambda$ in terms of the number of $a$'s,
$b$'s, $c$'s, and $d$'s in the diagram for $\lambda$ as follows:
\begin{eqnarray*}
\theta(\lambda) & = & \alpha(\lambda) - \beta(\lambda) +\gamma(\lambda) -
\delta(\lambda) \\
\theta(\lambda') & = & \alpha(\lambda) + \beta(\lambda) -\gamma(\lambda) -
\delta(\lambda) \\
|\lambda| & = & \alpha(\lambda) + \beta(\lambda) +\gamma(\lambda) +
\delta(\lambda).
\end{eqnarray*}
Thus we transform $\Phi(a,b,c,d)$ by sending $a \mapsto rsq$, $b \mapsto r^{-1}sq$,
$c \mapsto rs^{-1}q$, and $d \mapsto r^{-1}s^{-1}q$.  A straightforward computation
gives (\ref{And:result}). 

Our main result is a generalization of Theorem~\ref{main} and
Corollary~\ref{maincor}.  It is the corresponding product formula
in the case where we restrict the
the parts to some congruence class $(\mathrm{mod~}k)$ and we restrict the number of
times those parts can occur.   
Let $R$ be a subset of positive integers congruent to $i (\mathrm{mod~}k)$ and
let $\rho$ be a map from $R$ to the even positive integers.
Let $\mathrm{Par}(i,k;R, \rho)$ be the set of all partitions with parts
congruent to $i (\mathrm{mod~} k)$ such that if $r \in R$, then $r$ appears
as a part less than $\rho(r)$ times.  
Let $\Phi_{i,k;R,\rho}(a,b,c,d) = \sum_{\lambda} w(\lambda)$ where the sum is
over all partitions in $\mathrm{Par}(i,k;R,\rho)$.

For example, $\mathrm{Par}(1,1;\emptyset, \rho))$ 
is $\mathrm{Par}$, the set of all partitions.  Also, if we let $R$ be the set of
all positive integers and $\rho$ map every positive integer to $2$, then 
$\mathrm{Par}(1,1;R,\rho)$ is the set of all partitions with distinct parts.  These are the two
cases found in Theorem~\ref{main} and Corollary~\ref{maincor}.

\begin{theorem}
\label{theorem}
\[
\Phi_{i,k;R,\rho}(a,b,c,d) = ST
\]
where
\[
S = \prod_{j=1}^\infty 
\frac{(1+a^{\lceil \frac{(j+1)k+i}{2} \rceil}b^{\lfloor \frac{(j+1)k+i}{2} \rfloor}
c^{\lceil \frac{jk+i}{2} \rceil}d^{\lfloor \frac{jk+i}{2} \rfloor})
}
{(1-a^{\lceil \frac{jk+i}{2} \rceil}b^{\lfloor \frac{jk+i}{2} \rfloor}
c^{\lceil \frac{jk+i}{2} \rceil}d^{\lfloor \frac{jk+i}{2} \rfloor})
(1-a^{jk}b^{(j-1)k}c^{jk}d^{(j-1)k})}
\]
and
\[
T = \prod_{r \in R}(1-
a^{\lceil \frac{r}{2} \rceil \frac{\rho(r)}{2}}
b^{\lfloor \frac{r}{2} \rfloor \frac{\rho(r)}{2}}
c^{\lceil \frac{r}{2} \rceil \frac{\rho(r)}{2}}
d^{\lfloor \frac{r}{2} \rfloor \frac{\rho(r)}{2}})
\]
\end{theorem}

\section{Combinatorial Proof of these Results}

The proof of Theorem~\ref{theorem} is a slight modification of the proof of
Theorem~\ref{main} and Corollary~\ref{maincor}.  For clarity, we will first
give the argument in the special case where we consider all partitions and
partitions with distinct parts and then we will mention how  the proof can
be modified to work in general.

\begin{pf1}
\end{pf1}
Consider the following class of partitions:
\[
\mathcal{R} = \{ \lambda \in \mathrm{Par}: \lambda_{2i-1} - \lambda_{2i} \leq 1 \}.
\]
We are restricting the difference between
a part of $\lambda$ which is at an odd level and the following part of
$\lambda$ to be at most $1$.  

To find the generating function  for partitions in
$\mathcal{R}$ under weight $w(\lambda)$ we will decompose $\lambda \in
\mathcal{R}$ into blocks of height $2$,
$\{(\lambda_1,\lambda_2),(\lambda_3,\lambda_4), \dots \}$.  (In order to do this
if we have an odd
number of parts, add one
part equal to $0$.)
Since the difference of parts is
restricted  to either $0$ or $1$ 
at odd levels, we can only get two types of block.  For any $k \geq 1$, we can
have a block with two parts of length $k$, i.e. $(k,k)$. Call this Type I.
In addition, for any $k \geq 1$, we can have a
block with one part of length $k$ and then other of length $k-1$, i.e. $(k,k-1)$.  
Call this Type II.  

\begin{figure}
\centering
\includegraphics{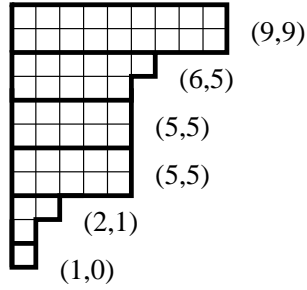}
\caption{$\lambda = (9,9,6,5,5,5,5,5,2,1,1)$ decomposes into blocks $\{(9,9),
(6,5), (5,5), (5,5), (2,1), (1,0)\}$}\label{fig1}
\end{figure}

In fact, partitions in $\mathcal{R}$ correspond uniquely to a
multiset of blocks of Type I and II with at most one block of Type II for each
length $k$. Figure~\ref{fig1} shows an example of such a decomposition.

To calculate the generation function for $\mathcal{R}$, it remains to calculate
the weigths of our blocks.  The blocks of Type I get filled as follows:
\[
\begin{array}{cccccccccccccccccc}
a&b&a&b& \dots &a&b & &\mathrm{or}& & a&b&a&b& \dots &a&b&a \\
c&d&c&d& \dots &c&d &  & &  & c&d&c&d& \dots &c&d&c
\end{array}
\]
depending on the length of the blocks.  Therefore they have weights
$a^{j}b^{j}c^{j}d^{j}$ or $a^{j}b^{j-1}c^{j}d^{j-1}$.

The blocks of Type II get filled as follows:
\[
\begin{array}{cccccccccccccccccc}
a&b&a&b& \dots &a&b&a & &\mathrm{or}& & a&b&a&b& \dots &a&b \\
c&d&c&d& \dots &c&d& &  & &  & c&d&c&d& \dots &c&
\end{array}
\]
depending on the length of the blocks.  Therefore they have weights
$a^{j}b^{j-1}c^{j-1}d^{j-1}$ or $a^{j}b^{j}c^{j}d^{j-1}$.

So we have the following generating function:
\[
\sum_{\lambda \in \mathcal{R}} w(\lambda) = 
\prod_{j=1}^\infty \frac{(1+a^{j}b^{j-1}c^{j-1}d^{j-1})(1+a^{j}b^{j}c^{j}d^{j-1})}
{(1-a^{j}b^{j}c^{j}d^{j})
(1-a^{j}b^{j-1}c^{j}d^{j-1})}.
\]

Notice that $\sum_{\lambda \in \mathcal{R}} w(\lambda)$ contains all the terms in
$\Phi(a,b,c,d)$ except for 
\[\prod_{j=1}^\infty \frac{1}
{1-a^{j}b^{j}c^{j-1}d^{j-1}}.\]  

Let $\mathcal{S}$ be the set of partitions whose conjugates have only odd
parts each of which is repeated an even number of times.  
We give a bijection $f:\mathrm{Par} \rightarrow \mathcal{R} \times \mathcal{S}$, such that $\mathcal{S}$ contributes exactly the
missing terms.  

\begin{figure}
\centering
\includegraphics{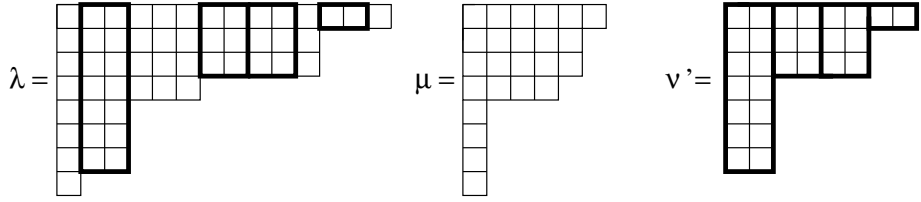}
\caption{$\lambda = (14,11,11,6,3,3,3,1)$ and $f(\lambda) = (\mu, \nu')$ where
$\nu = (7,7,3,3,3,3,1,1)$ and $\mu = (6,5,5,4,1,1,1,1)$}\label{fig2}
\end{figure}

Our bijection goes as follows, given a partition $\lambda$, let $f(\lambda) =
(\mu, \nu')$ where $\nu$ is the partition with $\lambda_{2i-1}-\lambda_{2i}$ parts
equal to $2i-1$ if $\lambda_{2i-1}-\lambda_{2i}$ is even and
$\lambda_{2i-1}-\lambda_{2i} - 1$ parts equal to $2i-1$ if
$\lambda_{2i-1}-\lambda_{2i}$ is odd, and where $\mu$ is defined by $\mu_i =
\lambda_i-\nu_i'$.   The function $f$ removes as many blocks of width 2 and odd
height as possible from $\lambda$. (Call these blocks of Type III.)  The remainder is a partition
$\mu$ such that the difference of parts at an odd level is either $0$ or $1$, that
is, a partition in $\mathcal{R}$.  The blocks of Type III which are removed are
joined together to give $\nu'$.  Clearly $f$ is a bijection.  
An example is shown in Figure~\ref{fig2}.

Now we examine the relationship between $w(\lambda)$, $w(\mu)$, and $w(\nu')$.
Consider the blocks of Type III in $\lambda$.  They always have weight
$a^jb^{j-1}c^jd^{j-1}$ regardless of whether their first column contains $a$'s and
$c$'s or $b$'s and $d$'s.  This is also the weight of the blocks when they are
placed in $\nu'$.  
Hence $w(\nu')$ is the product of the entries in the
diagram of $\lambda$ which are
removed to get $\mu$.

Moreover, since we are removing columns of width $2$, the entries in the squares of
the diagram of
$\lambda$ that correspond to squares in the diagram on $\mu$ do not change when
$\nu'$ is removed.  This implies that $w(\lambda) = w(\mu)w(\nu')$ and the result
follows.  
\qed

\begin{pf2}
\end{pf2}
To obtain this corollary, consider the following bijection.  Let $\mathcal{D}$
denote the set of partitions with distinct parts and let $\mathcal{E}$ denote the
set of partitions whose parts appear an even number of times.  Then we have a
bijection $g: \mathrm{Par} \rightarrow \mathcal{D} \times \mathcal{E}$ with 
$g(\lambda) = (\mu,\nu)$ defined as follows.
Suppose $\lambda$ has $k$ parts equal to $i$. 
If $k$ is even then $\nu$ has $k$ parts equal to $i$, and if $k$ is odd then
$\nu$ has $k-1$ parts equal to $i$.
The parts of $\lambda$ which were not removed to form $\nu$, at most one of each
cardinality, give $\mu$.
It is clear that under this bijection, $w(\lambda) = w(\mu)w(\nu)$.  

Now using the decompostion from the proof of Theorem~\ref{main}, partitions in 
$\mathcal{E}$ have a decomposition which only uses blocks of Type I.  
Hence we get that
\[
\Phi(a,b,c,d) = \Psi(a,b,c,d)\prod_{j=1}^\infty 
\frac{1}
{(1-a^{j}b^{j}c^{j}d^{j})
(1-a^{j}b^{j-1}c^{j}d^{j-1})}
\]
and the result follows.
\qed
	
\medskip	
The proof of our main result follows by the same argument with a modification to
the sizes of the blocks.

\begin{pf3}
\end{pf3}
First we find the generation function $S = \Phi_{i,k;\emptyset,\rho}(a,b,c,d)$ without
any restriction on the number of times each part may occur.  This is done by
using
Type I blocks with two parts each of length $jk+i$, for $j \geq 1$,  
Type II blocks with two
parts, one of length $jk+i$ and one of length $(j-1)k +i$, for $j \geq 1$, and  
Type III blocks which are rectangular with width $2k$ and odd height.

Next, we notice a bijection,
analogous to the one in the proof of Corollary~\ref{maincor}, between
$\mathrm{Par}(i,k;\emptyset,\rho)$ and 
$\mathrm{Par}(i,k;R,\rho) \times
\mathcal{T}$ where $\mathcal{T}$ is the set of all partitions with parts, $r \in
R$ and occuring a multiple of $\rho(r)$ times.  Since the generating function
for $\mathcal{T}$ is
\[
T^{-1} = \prod_{r \in R}\frac{1}{1-
a^{\lceil \frac{r}{2} \rceil \frac{\rho(r)}{2}}
b^{\lfloor \frac{r}{2} \rfloor \frac{\rho(r)}{2}}
c^{\lceil \frac{r}{2} \rceil \frac{\rho(r)}{2}}
d^{\lfloor \frac{r}{2} \rfloor \frac{\rho(r)}{2}}} 
\]
we get that
$S = \Phi_{i,k;R,\rho}(a,b,c,d)T^{-1}$ and the result follows.
\qed

\begin{ack}
\end{ack}
The author would like to thank Richard Stanley for encouraging her to work on
this problem.

\bibliography{abcd}

\end{document}